\documentclass[aps,prl,final,twocolumn,showpacs]{revtex4}

\usepackage{amsmath,amssymb}
\usepackage{graphicx}
\usepackage{epsf,url}

\newtheorem{theorem}{Theorem}
\newtheorem{prop}{Proposition}

\newtheorem{defn}{Definition}
\newtheorem{conjecture}{Conjecture}

\begin{document}

\title{Discrete tomography: Magic numbers for $N$-fold
  symmetry}

\author{Christian Huck, Markus Moll and Johan Nilsson}
\affiliation{Fakult\"{a}t f\"{u}r Mathematik, 
Universit\"{a}t Bielefeld,
Postfach 100131, 33501 Bielefeld, Germany}

\begin{abstract}
 We consider the problem of distinguishing 
  convex subsets of $n$-cyclotomic model sets $\varLambda$ by
  (discrete parallel) X-rays in prescribed
  $\varLambda$-directions. In this context, a `magic number'
  $m_{\varLambda}$ has the property that any two convex subsets of $\varLambda$  can be distinguished by
  their X-rays in any set of $m_{\varLambda}$ prescribed
  $\varLambda$-directions. Recent calculations suggest that (with one exception
  in the case $n=4$) the least possible
  magic number for $n$-cyclotomic model sets might just be $N+1$, where $N=\operatorname{lcm}(n,2)$.
 \end{abstract}


\maketitle

\section{Introduction}

{\em Discrete tomography} (DT) is concerned with the 
inverse problem of retrieving information about some {\em finite}
set in Euclidean space from 
(generally noisy) information about its slices. One important problem
is the {\em unique reconstruction} of a finite point set in Euclidean $3$-space 
from its {\em $($discrete parallel\/$)$ X-rays} in a small number of directions, where the
  {\em X-ray} of the finite set in a certain direction is the {\em line sum
  function} giving the
number of points in the set on each line parallel to this direction. 

The interest in the discrete tomography of planar Delone sets $\varLambda$
with long-range order is motivated by the requirement in materials
science for the unique reconstruction of solid state materials like {\em quasicrystals}
slice by slice from their images under quantitative {\em high
  resolution transmission electron microscopy} (HRTEM). In fact,
 a technique
is described in~\cite{ks,sk}, which can, for certain 
crystals, effectively measure the number of atoms lying on densely occupied 
columns. Since we aim at a slice by slice approach, it is
sufficient to study the DT of planar Delone sets with long-range
order.

Since different finite subsets of $\varLambda$ may have the same X-rays in several
directions, 
one is naturally interested in conditions
to be imposed
on the set of directions together with restrictions on
the possible finite subsets of
$\varLambda$ such that the latter phenomenon
cannot occur. Here, we consider the {\em convex subsets}
of so-called {\em strong cyclotomic Delone sets} $\varLambda$ and
present new results on the problem of
distinguishing convex subsets of  $\varLambda$ by X-rays in prescribed
$\varLambda$-directions. It is well known~\cite{HS} that there are four prescribed $\varLambda$-directions such that any two convex subsets of
$\varLambda$ can be distinguished by the corresponding X-rays,
whereas less than four $\varLambda$-directions never suffice for this
purpose. Moreover, there is a finite number
$m_{\varLambda}$ such that any two convex subsets of $\varLambda$
can be distinguished by their X-rays in {\em any} set of
$m_{\varLambda}$ prescribed
$\varLambda$-directions. It was essentially shown in~\cite{HS} that the
least possible `magic numbers' $m_{\varLambda}$ in the
case of the practically most relevant examples of strong $n$-cyclotomic Delone
sets $\varLambda$ with $n=5$, $8$ and $12$ only depend on $n$ and are (in that
 order) $11$, $9$ and $13$; see also~\cite{H} for a gentle
 introduction. This extended a well-known result of Gardner
and Gritzmann~\cite{GG} on the corresponding
problem for the crystallographic cases $n=3,4$ (with least possible number $m_{\varLambda}=7$ in
both cases) to
cases that are relevant in quasicrystallography. With the exception of the case 
$n=4$, the explicit calculations above indicate a very simple
relation between $n$ and the associated least possible `magic
numbers', i.e.\ for the cases $n=3,5,8$ and $12$ the magic number is
just $N+1$, where $N=\operatorname{lcm}(n,2)$. In this short note, we present further computational evidence that this might be
true in general.  For detailed proofs and an extensive list of references, we refer the reader
to~\cite{HS}.

\section{Cyclotomic Delone sets}

Throughout, the Euclidean plane is 
identified with the complex numbers. We denote by $K_{\varLambda}$
the smallest subfield of $\mathbb C$ that contains the rational
numbers as well as the union 
of $\varLambda-\varLambda$ 
and its image $\overline{\varLambda-\varLambda}$ under complex
conjugation. For $n\in \mathbb{N}$, we always let $\zeta_n = e^{2\pi
    i/n}$, a primitive $n$th root of unity in $\mathbb C$. Then, the smallest
  subfield of $\mathbb C$ that contains the rational
numbers as well as $\zeta_n$  is the
  $n$th cyclotomic field denoted by $\mathbb Q(\zeta_n)$. The latter
  is just the $\mathbb Q$-span of the $n$th roots of unity and thus
  contains the $\mathbb Z$-span $\mathbb Z[\zeta_n]$ of the $n$th
  roots of unity. Note that $\mathbb Q(\zeta_n)=\mathbb Q(\zeta_N)$,
  where $N=\operatorname{lcm}(n,2)=[n,2]$ is the least common multiple  of $n$ and
$2$. Recall
  that a {\em homothety}\/ of the complex plane  is given by $z \mapsto \lambda z + t$, where
$\lambda \in \mathbb R$ is positive and $t\in \mathbb C$. For our purposes, the following rather abstract definition is
  convenient. 

\begin{defn}\label{algdeldef}
Let $n\geq 3$. A Delone set $\varLambda\subset\mathbb C$ is called a
{\em strong $n$-cyclotomic
  Delone set} if it satisfies the following properties:

\begin{center}
\begin{tabular}{rl}
$($$n$-Cyc$^*$$)$&
$K_{\varLambda}= \mathbb Q(\zeta_n)$.\\
$($Hom\/$)$& \begin{minipage}[t]{6.5cm} For any finite subset $F$ of $K_{\varLambda}$, there is a
homothety $h$ of the complex plane that maps the
elements of $F$ to $\varLambda$.
\end{minipage}
\end{tabular}
\end{center}
Further,
$\varLambda$ is called a {\em strong cyclotomic Delone set} if it is a
strong 
$n$-cyclotomic Delone set for a suitable $n\geq 3$.
\end{defn} 
This includes many of the commonly used mathematical models of slices
that occur in quasicrystallography, i.e.\ the {\em $n$-cyclotomic model
  sets} $\varLambda\subset\mathbb Z[\zeta_n]$; cf.~\cite{BG,HS}. For suitable choices of
the window, these sets have $N$-fold rotational symmetry. They range from periodic examples like
the fourfold square lattice ($n=4$) or the sixfold triangular lattice
($n=3$) to nonperiodic examples like the vertex set of the tenfold T\"ubingen triangle
tiling ($n=5$), the eightfold Ammann-Beenker tiling ($n=8$) or the twelvefold
shield tiling ($n=12$);
see Figure~\ref{fig:tilingupolygon}. The vertex sets of Penrose
tilings fail to be cyclotomic model sets but are still strong $5$-cyclotomic Delone sets~\cite{bh}. 

\section{Discrete tomography}

The
unit circle in $\mathbb C$ is
denoted by $\mathbb{S}^{1}$ and its elements are also called
{\em directions}. 

\begin{defn}
Let $F$ be a finite subset of\/ $\mathbb C$, let $u\in
\mathbb{S}^{1}$ be a direction, and let $\mathcal{L}_{u}$ be the set
of lines in the complex plane in direction $u$. Then the {\em
  (discrete parallel)}\/ {\em X-ray} of $F$ {\em in direction} $u$ is
the function $X_{u}F: \mathcal{L}_{u} \rightarrow
\mathbb{N}_{0}=\mathbb{N} \cup\{0\}$, defined by $$X_{u}F(\ell) =
\operatorname{card}(F \cap \ell\,).$$
Furthermore, we say that the elements
of a collection $\mathcal{F}$ of finite subsets of\/
$\mathbb C$ are {\em determined} by the X-rays in the directions
of a finite set $U\subset\mathbb{S}^{1}$ of directions if different elements of $\mathcal{F}$ cannot have the same
X-rays in the directions of $U$. 
\end{defn}

Let $\varLambda\subset\mathbb C$ be a fixed Delone set. Obviously, only {\em
  $\varLambda$-directions} (directions
parallel to nonzero elements of the {\em difference set}\/ 
$\varLambda-\varLambda$) are reasonable. One can see that, for a
strong $n$-cyclotomic Delone set $\varLambda$, the set of
$\varLambda$-directions coincides with the set of $\mathbb
Q(\zeta_n)$-directions (defined analogously). Note that the present text is only concerned with
the {\em uniqueness problem} of determining the elements of a large collection of finite
subsets of $\varLambda$ by few X-rays in prescribed
$\varLambda$-directions. For the algorithmic 
{\em reconstruction problem} in the quasicrystallographic setting,
see~\cite{H2}. 

One can easily see that no finite set of pairwise nonparallel
$\varLambda$-directions suffices in order to determine the whole class
of finite subsets of a (strong) cyclotomic Delone set $\varLambda$
(in contrast to strong cyclotomic Delone sets, cyclotomic Delone sets
only satisfy the weaker condition $K_{\varLambda}\subset\mathbb Q(\zeta_n)$)  by the corresponding
X-rays. It has proven to be most fruitful to focus on the {\em convex
  subsets}\/ of (strong) cyclotomic Delone
sets. The latter are bounded (and thus finite) subsets $C$ of $\varLambda$ with $C =
\operatorname{conv}(C)\cap \varLambda$, where $\operatorname{conv}(C)$
denotes the convex hull of $C$.

\section{Determining convex sets by X-rays}

\begin{defn}\label{defupol}
For  a finite set $U\subset \mathbb{S}^{1}$ of
directions, a nondegenerate convex polygon $P\subset\mathbb C$ is
called a {\em $U$-polygon} if it has the property that whenever $v$ is
a vertex of $P$ and $u\in U$, the line in the complex plane in direction $u$
which passes through $v$ also meets another vertex $v'$ of $P$. For a
subset $S\subset \mathbb C$, $P$ is
called a 
{\em $U$-polygon in $S$}, if its vertices lie in $S$.
\end{defn}

One has the following fundamental 
result. 

\begin{center}
\begin{figure}
\centerline{\epsfysize=0.48\textwidth\epsfbox{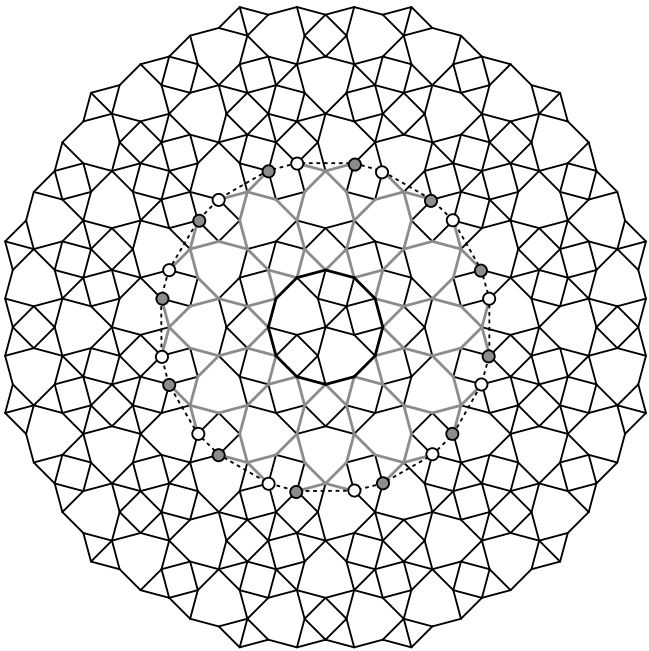}}
\caption{The construction of a $U$-polygon in
  the vertex set $\varLambda$ of the twelvefold 
  shield tiling (as described in the paragraph preceding Proposition~\ref{prop1}), where $U$ is the set of
  $N=12$ $\varLambda$-directions
 given by the edges and diagonals of the central regular
 dodecagon. Here, the points obtained by translating
 the regular $N$-gon are already vertices of the tiling and not only
 elements of $\mathbb
 Q(\zeta_{n})$. The vertices of $\varLambda$ in the interior of
 the  $U$-polygon together with the vertices indicated by the white and grey dots, respectively,
give two different convex subsets of $\varLambda$ with the same X-rays in the
directions of $U$.}
\label{fig:tilingupolygon}
\end{figure}
\end{center}

\begin{theorem}{\rm \cite[Fact 5.3]{HS}}\label{characungen}
Let $\varLambda$ be a (strong) cyclotomic Delone set and let $U\subset\mathbb{S}^{1}$ be a set of two or more pairwise nonparallel $\varLambda$-directions. The following statements are equivalent:
\begin{itemize}
\item[\rm (i)]
The convex subsets of $\varLambda$ are determined by the X-rays in the directions of $U$.
\item[\rm (ii)]
There is no $U$-polygon in $\varLambda$.
\end{itemize}
If, in addition, $\operatorname{card}(U)<4$, there is a $U$-polygon in $\varLambda$. 
\end{theorem}
Note
that the proof of direction
(ii)$\Rightarrow$(i) needs property
(Hom);  see Figure~\ref{fig:tilingupolygon} for an
illustration of the other (easy) direction
(i)$\Rightarrow$(ii). Thus, one is immediately lead to the investigation of $U$-polygons.

\section{Existence of $U$-polygons}

For $n\ge 3$, consider the regular $N$-gon $R$ inscribed in the unit circle, with one
vertex at $1$. Now attach $N$ translates of $R$ edge-to-edge  to $R$
in the obvious way. Then, one can easily verify that the convex hull $P$ of the resulting point
set is a $U$-polygon in $\mathbb Q(\zeta_n)$, where
$U$ is the set of $N$ pairwise nonparallel $\mathbb Q(\zeta_n)$-directions
 given by the edges and diagonals of $R$. Now let $\varLambda$ be a
 strong $n$-cyclotomic Delone set. Employing properties
 ($n$-Cyc$^*$) and (Hom), one obtains a $U$-polygon in
 $\varLambda$ for a set $U$ consisting of $N$
 $\varLambda$-directions;
compare Figure~\ref{fig:tilingupolygon}. In view of Theorem~\ref{characungen}, this proves that the magic number for strong
 $n$-cyclotomic Delone sets is at least $N+1$.   

\begin{prop}\label{prop1}
Let $\varLambda$ be a strong $n$-cyclotomic Delone set. Then there is a
$U$-polygon in $\varLambda$ with $\operatorname{card}(U)=N$. 
\end{prop}

\begin{defn}
Let $(t_1,t_2,t_3,t_4)$ be an ordered tuple of four distinct
elements of $\mathbb{R}\cup\{\infty\}$. Then, its {\em cross ratio}\/
$\langle t_1,t_2,t_3,t_4\rangle$ is %the nonzero real number defined by
$$
\langle t_1,t_2,t_3,t_4\rangle = \frac{(t_3 - t_1)(t_4 - t_2)}{(t_3 -
  t_2)(t_4 - t_1)}\in\mathbb R^*,
$$
with the usual conventions if one of the $t_i$ equals
$\infty$. 
\end{defn}
By construction, the
cross ratio of slopes of four pairwise nonparallel
$\varLambda$-directions is an element of the field
$K_{\varLambda}\cap\mathbb R$. In case of (strong) $n$-cyclotomic Delone
sets $\varLambda$, these cross ratios are thus elements of the field
$\mathbb Q(\zeta_n)\cap\mathbb R$.

Employing a blend of methods from the theory of
cyclotomic fields and previous results obtained by Gardner and Gritzmann~\cite{GG}, one obtains the
following deep result on $U$-polygons. The novelty here is the
assertion on the formula for the least possible upper bounds
$b_n$ and the largest possible sets $U$ of directions. These were found
by direct computation (the proof of~\cite[Thm.\ 5.6]{HS} immediately leads to an algorithm). 

\begin{theorem}{\rm \cite[Thm. 5.7]{HS}}\label{finitesetncr0gen}
Let $n\geq 3$ and let $\varLambda$ be a strong $n$-cyclotomic Delone set. Further, let $U\subset\mathbb{S}^1$ be a set of four or more pairwise
nonparallel $\varLambda$-directions and suppose the existence of a $U$-polygon. Then, the
cross ratio of slopes of any four directions of $U$, arranged in order
of increasing angle with the positive real axis, is an element of the
finite set $\mathcal{C}_{[2n,12]}(\mathbb Q(\zeta_n)\cap\mathbb R)
$ of numbers in the field $\mathbb Q(\zeta_n)\cap\mathbb R$ that can be written
in the form
$$
\frac{\Big(1-\zeta_{[2n,12]}^{k_1}\Big)\Big(1-\zeta_{[2n,12]}^{k_2}\Big)}{\Big(1-\zeta_{[2n,12]}^{k_3}\Big)\Big(1-\zeta_{[2n,12]}^{k_4}\Big)},
$$
where $(k_1,k_2,k_3,k_4)\in \mathbb{N}^4$ satisfies
$$
 k_3<k_1\leq
k_2<k_4\leq [2n,12]-1 \mbox{ and } k_1+k_2=k_3+k_4
.
$$
Moreover, $\operatorname{card}(U)$ is bounded above by a finite number
$b_n\in\mathbb N$ that only depends on $n$. With the exception of $b_4=7$, for $3\le
n\le 50$ and $n=61$, $b_n=N=[n,2]$ is best possible. Moreover, for
these values of $n$, if there is a $U$-polygon with
$\operatorname{card}(U)=N$, then there is a linear automorphism of the
complex plane that takes the directions from $U$ to a set of vectors
which when normalised are given by $e^{h\pi i/N}$, where
$h\in\{0,1,2,\dots,N-1\}$ (compare Figure~{\rm \ref{fig:tilingupolygon}}). 
\end{theorem}

\section{Result}

Theorems~\ref{characungen} and~\ref{finitesetncr0gen} now immediately imply our main result on the
determination of convex subsets of strong cyclotomic Delone sets.

\begin{theorem}{\rm \cite[Thm. 5.11]{HS}}\label{dtmain2}
Let $n\geq 3$ and let $\varLambda$ be a strong $n$-cyclotomic Delone set.
\begin{itemize}
\item[\rm (a)]
There are sets of four pairwise nonparallel $\varLambda$-directions
such that the convex subsets of $\varLambda$ are determined by the corresponding
X-rays. In addition, less than four pairwise nonparallel $\varLambda$-directions never
suffice for this purpose.
\item[\rm (b)]
There is a finite number $m_n\in\mathbb N$ that only depends on $n$ such that the convex subsets of
$\varLambda$ are determined by the X-rays in any set of $m_n$
pairwise nonparallel $\varLambda$-directions. With the exception of\/
$m_4=7$, for $3\le
n\le 50$ and $n=61$, $m_n=N+1$ is best possible. 
\end{itemize}
\end{theorem}

\section{Recipe and conjecture}

By Theorems~\ref{characungen} and~\ref{finitesetncr0gen} above, it
suffices for Theorem~\ref{dtmain2}(a) to take any set of four pairwise
nonparallel $\varLambda$-directions such that the cross ratio of their
slopes, arranged in order
of increasing angle with the positive real axis, is {\em not} an element of
the finite set $\mathcal{C}_{[2n,12]}\big(\mathbb
Q(\zeta_n)\cap\mathbb R\big)$ which can be computed directly. E.g., for the
vertex set of the shield tiling (a $12$-cyclotomic model set) below,
the set $\mathcal{C}_{24}\big(\mathbb Q(\sqrt{3})\big)$ of cross ratios to be avoided
is given by (note that the number $\frac{2}{\sqrt{3}}$ is mistakenly 
  missing in the corresponding list appearing in~\cite[Cor.\
  4.10(c)]{HS}) 

\begin{eqnarray*} 
&&
\Big\{8-4\sqrt{3},\frac{3+2\sqrt{3}}{6},\frac{-3+3\sqrt{3}}{2},\frac{2}{\sqrt{3}},\frac{3+\sqrt{3}}{4},\frac{2+\sqrt{3}}{3},\\
&&\hphantom{\Big\{}
3-\sqrt{3},\frac{4}{3},\frac{1+\sqrt{3}}{2},-2+2\sqrt{3},\frac{3}{2},\frac{3+\sqrt{3}}{3},\sqrt{3},\\
&&\hphantom{\Big\{}\frac{2+\sqrt{3}}{2},2,\frac{3+2\sqrt{3}}{3},\frac{3+\sqrt{3}}{2},1+\sqrt{3},3,\frac{6+2\sqrt{3}}{3},\\
&&\hphantom{\Big\{}2+\sqrt{3},4,3+\sqrt{3},\frac{5+3\sqrt{3}}{2}, 3+2\sqrt{3},4+2\sqrt{3},\\
&&\hphantom{\Big\{}6+3\sqrt{3},7+4\sqrt{3},8+4\sqrt{3}\Big\}.
\end{eqnarray*}
See~\cite[Cor. 4.10]{HS} for concrete results in the other practically important cases $n=5,8$
of
quasicrystallography.

\begin{conjecture}
With the exception of\/
$m_4=7$, the trivial lower bound 
$
m_n=N+1
$
is the magic number for strong $n$-cyclotomic Delone sets for all
$n\ge 3$, where $N=[n,2]$.
\end{conjecture}
The above exception can be explained as follows. For the two
crystallographic cases $n=3,4$, the corresponding $n$-cyclotomic model
sets, i.e.\ translates of the  triangular ($N=6$) resp.\ square
lattice ($N=4$), are affinely equivalent. It is thus clear that, as it is the case for the
triangular lattice, also the square lattice contains a $U$-polygon for
a set $U$ consisting of \emph{six}\/ (and not only $N=4$) lattice
directions; cf.~\cite{GG}. 

\section{Seven a universal magic number in 3D?}

It is certainly of interest to abandon the slice-by-slice
approach to the problem of uniquely reconstructing convex
subsets of 3D Delone sets by X-rays and
to turn to sets $U\subset\mathbb{S}^2$ of $\varLambda$-directions in
{\em general position}\/ (i.e.\ any three directions from $U$ span the
whole $3$-space) instead; compare the
approach to 3D reconstruction of atomic arrangements presented in~\cite{Saitoh}. The
following conjecture goes back to a continuous version formulated by
Gardner in 1995~\cite[Problem 2.1]{G}. 

\begin{center}
    \begin{figure}
        \begin{center}
           \includegraphics[scale=0.75]{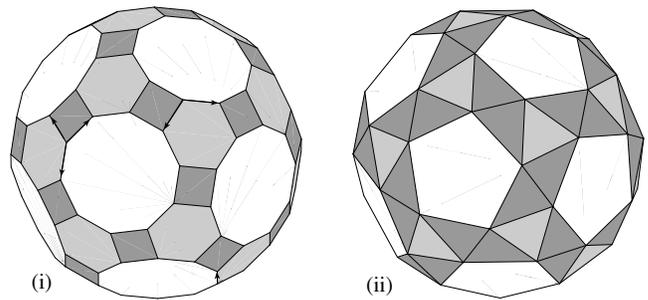}
        \end{center}
        \caption{The great rhombicosidodecahedron (i) (inscribed in
          the unit sphere) is a
        $U$-polyhedron (defined just as $U$-polygons) with respect to the indicated set $U$ of
        six directions in general position. Again, its vertices can be
        coloured with two colours, say white and grey, such that every
        white vertex
        corresponds to a grey vertex, whatever direction from
        $U$ is used. The convex hull of one colour class is shown in
        (ii). For symmetry reasons the vertices of (ii) and the
        vertices of its counterpart (the convex hull of the other
        colour class) have the same X-rays in the
        directions of $U$. Figure from~\cite[p.\ 66]{G}
        {\copyright}~Cambridge University Press.}
        \label{fig:lowerbound}
    \end{figure}
\end{center}

\begin{conjecture}\label{conj:gardner}
   Let $\varLambda\subset \mathbb R^3$ be a Delone set. Then, the
   convex subsets of  $\varLambda$ can be distinguished by
    the X-rays in any set of seven $\varLambda$-directions in general
    position.
\end{conjecture}
Figure~\ref{fig:lowerbound} immediately yields that the number seven in
Conjecture~\ref{conj:gardner} cannot be lowered to six.

\section*{Acknowledgements}

The authors would like to thank Franz G\"ahler for the help and
assistance with the computer calculations. This work was supported by the German Research Foundation (DFG) within
the CRC~701.

\end{document}